\documentclass[11pt]{article}

\batchmode

\usepackage{amsfonts}

\newtheorem{theorem}{Theorem}
\newtheorem{lemma}{Lemma}

\newtheorem{definition}{Definition}

\newtheorem{corr}{Corollary}

\def\lor{\mathop{\mbox{\Large $($}}}  
\def\ler{\mathop{\mbox{\Large $)$}}}  
\def\={\mathop{\;=\;}}  

\def\Largabs{\mathop{\mbox{\Large $|$}}}

\newcounter{aid}

\begin{document}

\title{\bf de Finetti Style Theorems With Applications to Network Analysis}

\author{ \vspace*{1mm} Andr\'as  Farag\'o \\
Department of Computer  Science      \\
The  University   of  Texas   at  Dallas\\
   Richardson,  Texas \\
   {\tt farago@utdallas.edu} }

\date{}
\maketitle

\begin{abstract} A classic and fundamental result about the decomposition of random sequences into a mixture of simpler ones is de Finetti's Theorem.
In its original form it applies to infinite  0--1 valued exchangeable sequences. 
Later it was extended and generalized in numerous directions. After reviewing this line of development, we present our new decomposition theorem, 
covering cases that have not been previously considered. 
We also introduce a novel way of applying these types of results in the analysis of random networks.
For self-containment, we provide the introductory exposition in more details than usual, with the intent of making it also accessible to readers 
who may not be closely familiar with the subject.
 \end{abstract}

\thispagestyle{empty}

\section{Introduction and Background}
\label{intro}

It has been a long standing pursuit in probability theory and its applications to express a random sequence as a mixture of simpler random sequences. 
The mixing is meant here in the probabilistic sense, that is, we select one among  the component sequences via some
probability distribution that governs the mixing, and then output the selected sequence in its entirety.
Equivalently, the distribution of the resulting sequence (i.e., the joint distribution of its entries) is a convex combination 
of the distributions of the component sequences. The distribution used for the selection  is often referred to as the {\em mixing measure}.

{\em Note:}  when we only want to represent a single random variable as a mixture, it is a much simpler case, discussed in
the well established statistical field  of {\em mixture 
models,} see Lindsay \cite{lindsay}. Here we are interested, however, in expressing random {\em sequences,} rather than just single random variables.

Which simple sequences can serve best as the components of mixing? Arguably, the simplest possible probabilistic structure that a random sequence can have is being a sequence of {\em independent, identically distributed (i.i.d.)} random variables. 
The mixture of such i.i.d.\ sequences, however, does not have to remain i.i.d. For example, the identically 0 and identically 1 sequences are both i.i.d., but if we mix them by selecting one of them with probability 1/2, then we get a sequence in which each term is either 0 or 1 with probability 1/2, but all of them are equal, so the entries are clearly not independent.

Since the joint distribution of any i.i.d.\ sequence is invariant to reordering the terms by any {\em fixed} permutation, therefore, the mixture must also
behave this way. The reason is that it does not matter whether we first apply a permutation $\sigma$ to each sequence and then select one of them, or first make the selection and apply the permutation afterwards to the selected sequence. The sequences with the property that their joint distribution is invariant to permutations are called {\em exchangeable:}

\begin{definition} \label{def1} {\bf (Exchangeable sequence)} A finite sequence $\xi=(\xi_1,\ldots,\xi_n)$ of random variables is called {\em exchangeable} 
if its joint distribution is invariant with respect to permutations. That is,
for any permutation 
$\sigma$ of $\{1,\ldots,n\}$, the joint distribution of $(\xi_{\sigma(1)},\ldots, \xi_{\sigma(n)} )$ is the same as the joint distribution of 
$(\xi_1,\ldots,\xi_n)$. An infinite sequence is called exchangeable if every finite initial segment of the sequence is exchangeable. 
\end{definition}
It means, an exchangeable sequence is stochastically indistinguishable from any permutation of itself. 
An equivalent definition is that if we pick $k$ entries of the sequence, then
their joint distribution depends only on $k$, but not on that which $k$ entries are selected and in which order. 
This also implies (with $k=1$) that each individual entry has the
same distribution. Sampling tasks often produce exchangeable sequences, since in most cases the order of the samples does not matter.

As a special case, the definition is satisfied by 
i.i.d.\ random variables, but not every exchangeable sequence is i.i.d. 
There are many examples to demonstrate this; a simple one can be obtained 
from geometric considerations:

{\sc Example 1.} Take a square in the plane, and divide it into two triangles by one of its diagonals. Select  one of the triangles 
at random with probability 1/2, and then pick $n$ uniformly random points from the selected triangle. These random points constitute  an exchangeable 
sequence, since their joint probability distribution remains the same, regardless of the order they have been produced. Furthermore, each individual point is 
uniformly distributed over the whole square, because it is uniformly distributed over a triangle, which is selected with equal probability 
from among the two triangles. On the other hand, the random points  are not independent, 
since if we know that a point falls in the interior of a given one of the two triangles, the all the others must fall in the same triangle. 

As we have argued before Definition~\ref{def1}, the mixing of i.i.d.\ sequences produces exchangeable sequences. 
A classical theorem of Bruno de Finetti, originally published in Italian \cite{definetti} in 1931, says
that the converse is also true for infinite binary sequences: every infinite exchangeable sequence of binary random variables can be
represented as a mixture of i.i.d.\ Bernoulli random variables (for short, a Bernoulli i.i.d.-mix). 
The result can be formally stated in several different ways, here is a frequently used one,
which captures the distribution of the exchangeable sequence as a mixture of binomial distributions:

\begin{theorem} \label{thm1} {\bf (de Finetti's Theorem -- distributional form)}
Let $X_1,X_2,\ldots$ be an infinite sequence of $\{0,1\}$-valued exchangeable random variables. Then there exists a probability measure $\mu$ 
(called mixing measure)
on $[0,1]$, such that for every positive integer $n$ and for any $x_1,\ldots,x_n \in \{0,1\}$ the following holds:
\begin{equation}\label{mix}
\Pr(X_1=x_1,\ldots,X_n=x_n) \= \int_0^1 p^s(1-p)^{n-s} d\mu(p)
\end{equation}
where $s=\sum_{i=1}^n x_i$. Furthermore, the measure $\mu$  is uniquely determined.
\end{theorem}

Note that the reason for using Stieltjes integral on the right-hand side of (\ref{mix}) is just to express discrete, continuous and mixed distributions 
in a unified format. For example, if the mixing measure $\mu$ is discrete, taking values $a_1,a_2,\ldots$ with probabilities $p_1,p_2,\ldots$, respectively,
then the integral becomes the sum
$
\sum_i a_i^s(1-a_i)^{n-s} p_i.
$
If the mixing measure is continuous and has a density function $\mu'$,  then the integral  becomes the ordinary integral 
$\int_0^1 p^s(1-p)^{n-s} \mu'(p) dp$. The Stieltjes integral expression contains all these special cases in a unified format, 
including mixed distributions, as well. 

Another often seen form of the theorem emphasizes that $X_1,X_2,\ldots$ becomes an i.i.d.\  Bernoulli sequence, whenever we {\em condition} 
on the value $p=\Pr(X_i=1)$, as presented below: 

\begin{theorem} \label{thm1b} {\bf (de Finetti's Theorem -- conditional independence form)}
Let $X_1,X_2,\ldots$ be an infinite sequence of $\{0,1\}$-valued exchangeable random variables. Then there exists a random variable $\eta$,
taking values in $[0,1]$, such that for every $p\in [0,1]$,
for every positive integer $n$ and for any $x_1,\ldots,x_n \in \{0,1\}$ the following holds:
\begin{equation}\label{mixb}
\Pr(X_1=x_1,\ldots,X_n=x_n\;|\;\eta=p) \= p^s(1-p)^{n-s} 
\end{equation}
where $s=\sum_{i=1}^n x_i$. Furthermore, $\eta$ is the limiting fraction of the number of ones in the sequence (the empirical distribution):
$$
\eta \= \lim_{n\rightarrow\infty} \frac{1}{n}\sum_{i=1}^n X_i
$$
\end{theorem}

It is interesting that the requirement of having an {\em infinite} sequence is essential; for the finite case counterexamples are known, 
see, e.g., Stoyanov  \cite{stoyanov}.
(Note that even though equations (\ref{mix}) and (\ref{mixb}) use a fixed finite $n$, the theorem requires it to hold for every $n$.)
On the other hand, approximate versions exist for finite sequences, see Section~\ref{finite}. It is also worth noting that the proof is far from easy.
An elementary proof was published by Kirsch \cite{kirsch} in 2019, but this happened 88 years after the original paper. 

\subsection{Philosophical Interpretation of de Finetti's Theorem}
\label{phyl}

The concept of probability has several philosophical interpretations (for a survey, see \cite{stanford}). An appealing aspect of de Finetti's Theorem that it builds a bridge between two major conflicting interpretations: the {\em frequentist} and the {\em subjective} interpretations. (The latter is also known as {\em Bayesian} interpretation.)  Let us briefly explain these trough the simple experiment of coin flipping.

The frequentist interpretation of probability says that there exists a real number $p\in [0,1]$, such that if we keep flipping 
the same coin independently, then the relative frequency of heads converges to $p$, and this value gives us the probability of heads. 
In this sense the probability is an {\em objective} quantity, even when we may not know its exact value. 
Most researchers accept this interpretation, since it is in good agreement with experiments, and provides a common-sense, testable concept. 
In some cases, however, it does not work so well, such as when we deal with a one-time event which cannot be indefinitely repeated.
For example, it is hard to assign a precise meaning to a statement like ``candidate X will win the election tomorrow with probability 52\%."

In contrast, the subjective (Bayesian) interpretation denies the objective existence of probability. Rather, it says that the concept 
only expresses one's subjective expectation that a certain event happens. For example, there is no reason to {\em a priori} assume that if 
among the first 100 coin flips we observed, say, 53 heads, then similar behavior has to be expected among the next 100 flips. 
If we still assume that the order in which the coin flips are recorded does not matter, then what we see is just an exchangeable sequence of binary values,
but possibly no convergence to a {\em constant.}

Which interpretation is right? The one that de Finetti favored (see \cite{definetti2}), against the majority view, was the subjective interpretation. Nevertheless, his theorem provides 
a nice bridge between the two interpretations, in the following way. Consider two experiments:

(1) {\em Bayesian:} Just keep flipping a coin and record the results. Do not presuppose the existence of a  probability to which the 
relative frequency of heads converges, but still assume that the order of recording does not matter. 
Then what we obtain is an exchangeable sequence, but no specific objective probability. 

(2) {\em Frequentist:} 
Assume that an objective probability $p$ of heads does exist, but we do not know its value exactly, so we consider it as a random quantity, drawn from some 
probability distribution $\mu$. Then the experiment will be this: draw $p$ from the distribution $\mu$, fix it, and then keep flipping a coin that
has probability $p$ of heads, on the basis that this probability $p$ objectively exists.

Now de Finetti's Theorem states that the results of the above two experiments are {\em indistinguishable:} an exchangeable sequence of coin flips cannot 
be distinguished from  a mix of Bernoulli sequences. In this sense, the conflicting interpretations do not lead to conflicting experimental results,
so the theorem indeed builds a bridge between the subjective and frequentist views. This is a reassuring reconciliation between the conflicting 
interpretations! 

We need to note, however, that the above argument is only guaranteed to work if the sequence of coin flips is infinite. As already mentioned earlier, 
for the finite case Theorem~\ref{thm1} does not always hold. This anomaly with finite sequences may be explained with the fact that 
the frequentist probability, as the {\em limiting} value of the relative frequency, is only meaningful if we can consider infinite sequences.

\section{Generalizations/Modifications of de Finetti's Theorem}

As the original theorem was published almost a century ago, and has been regarded a fundamental result since then, it is not surprising that numerous extensions, generalizations and modifications were obtained  over the decades. Below we briefly survey some of the typical clusters of the development.

\subsection{Extending the Result to More General Random Variables}

The original theorem, published in 1931, refers to binary random variables. In 1937, de Finetti himself showed \cite{definetti2} that it also holds for real valued random variables. This was extended to much more general cases in 1955 by Hewitt and Savage \cite{hewitt}. They allow random variables that 
take values from  a variety of very general spaces; one of the most general examples is a Borel measurable space (Borel space, for short; 
see the definition and explanation of related concepts in Appendix A). This space includes all cases that are likely to be encountered in applications.

To formally present the generalization of de Finetti's Theorem in a form similar to Theorem~\ref{thm1}, let $S$ denote the space from which the random variables take their values, and let ${\cal P}(S)$ be the family of all probability distributions on $S$.

\begin{theorem} \label{thm2} {\bf (Hewitt-Savage Theorem)}
Let $X_1,X_2,\ldots$ be an infinite sequence of $S$-valued exchangeable random variables, where $S$ is a Borel measurable space. 
Then there exists a probability measure $\mu$ on ${\cal P}(S)$, such that for every 
positive integer $n$ and for any measurable $A_1,\ldots,A_n \subseteq S$ the following holds:
$$
\Pr(X_1\in A_1,\ldots,X_n\in A_n) \= \int_{{\cal P}(S)} \pi(A_1)\cdot\ldots\cdot \pi(A_n)\, d\mu(\pi)
$$
where $\pi\in {\cal P}(S)$ denotes a random probability distribution from ${\cal P}(S)$, drawn according to $\mu$. 
Furthermore, the mixing measure $\mu$ is uniquely determined.
\end{theorem}
Less formally, we can state it this way: an infinite sequence of $S$-valued exchangeable random variables is an $S$-valued i.i.d.-mix, whenever 
$S$ is a Borel measurable space. 
Note that here $\mu$ selects a random distribution from ${\cal P}(S)$, which may be a complex object, while in Theorem~\ref{thm1} this random
distribution is determined by a single real parameter $p\in [0,1]$.

An interesting related result (which was actually published in the same paper \cite{hewitt}) is called {\em Hewitt-Savage 0--1 Law.}
Let $X=(X_1,X_2,\ldots)$ be an infinite i.i.d.\ sequence. Further, let $\cal E$ be an event that is determined by $X$. We say that $\cal E$ is {\em symmetric}
(with respect to $X$)
if the occurrence or non-occurrence of $\cal E$ is not influenced by permuting any finite initial segment of $X$. For example, the event that ``$X$ 
falls in a given set $A$ infinitely many times" is clearly symmetric, as it is not influenced by permuting any finite initial segment of $X$. The
Hewitt-Savage 0--1 Law says that any such symmetric event has probability either 0 or 1. 


As an illustration for Theorem~\ref{thm2},  consider the following example:

\noindent
{\sc Example 2.} Assume we put a number of balls into an urn. Each ball has a color, one of $t$ possible colors (the number of colors may be infinite). 
Let $k_i$ be the initial number of 
balls of color $i$ in the urn,  where the $k_i$ are arbitrary fixed non-negative integers. Consider now the following process: draw a ball randomly from 
the urn, let is color be denoted by $X_1$. Then put back {\em two} balls of the same color $X_1$ in the urn. Keep repeating this experiment by 
always drawing a ball randomly from the urn, and each time putting back two balls of the same color as that of the currently drawn ball. 
Let $X_1, X_2, X_3, \ldots$ denote the 
random sequence of obtained colors. This is called a $t$-color P\'olya urn scheme and it is known that the generated sequence 
$X_1, X_2, X_3, \ldots$ is exchangeable, see Hill, Lane and Sudderth \cite{hill}. Then, by Theorem~\ref{thm2}, the sequence can be represented as an i.i.d.-mix.
Note that just from the definition of the urn process this fact may be far from obvious. 

In view of the generality of Theorem~\ref{thm2}, one may push further: does the result hold for {\em completely arbitrary} random variables? After all, 
it does not seem self-explanatory why they need to take their values from a Borel measurable space. The most general target space that is allowed for 
random variables is a general {\em measurable space}, see the definition in Appendix A. One may ask: does the theorem remain true for random variables that 
take their values from {\em any} measurable space?

Interestingly, the answer is no. Dubins and Freedman \cite{dubins} prove that Theorem~\ref{thm2} does not remain true 
for this completely general case, so some structural restrictions are indeed needed,
although these restrictions are highly unlikely to hinder any application. A challenging question, however, still remains:  how can we explain the need for
such restrictions in the context of the philosophical interpretation outlined in Section~\ref{phyl}? Let us just mention, without elaborating on details, that restricting the general measurable space to a Borel measurable space means a topological 
restriction\footnote{For background on topological spaces we refer to the literature, see, e.g., Willard \cite{willard}.}. 
At the same time, topology can be viewed as a (strongly)  abstracted version
of geometry. In this sense, we can say that de Finetti-style theorems require that, no matter how remotely, we still have to somehow relate to
 the real world: at least some very abstract version of geometry is indispensable.

\subsection{Modifying the Exchangeability Requirement}
\label{modifying}

There are numerous results that prove some variant of  de Finetti's theorem (and of its more general version, the Hewitt-Savage theorem)
for random structures that satisfy some symmetry 
requirement similar to exchangeability. For a survey see Aldous \cite{aldous} and  Kallenberg \cite{kallenberg}. Here we present two characteristic examples.

\medskip
\noindent
{\bf Partially exhangeable arrays.} Let $X=[X_{i,j}], \, 1\leq i,j<\infty,$ be a doubly infinite array (infinite matrix) of random variables, taking values from a Borel 
measurable space $S$. Let $R_i, C_j$ denote the
$i^{th}$ row and $j^{th}$ column of $X$, respectively. We say that $X$ is {\em row-exchangeable,} if the sequence $R_1,R_2,\ldots$ is exchangeable. Similarly, 
$X$ is {\em column-exchangeable} if $C_1,C_2,\ldots$ is exchangeable. Finally, $X$ is {\em row and column exchangeable (RCE),} if $X$ is both row-exchangeable
and column-exchangeable. Observe that RCE is a weaker requirement than demanding that all entries of $X$, listed as a single sequence, form an 
exchangeable sequence. For RCE arrays Aldous \cite{aldous2} proved a characterization, which contains the de Finetti (in fact, the Hewitt-Savage) theorem
as a special case. We use the notation $X=_d Y$ to express that the random variables $X,Y$ have the same distribution.

\begin{theorem} {\bf (Row and column exchangeable (RCE) arrays)} If $X$ is an RCE array, then there exists independent random variables 
$\alpha, \xi_i,\eta_j, \lambda_{i,j}, \, 1\leq i,j<\infty,$ such that all of them are uniformly distributed on $[0,1]$, and there exists a 
measurable function\footnote{See the definition in Appendix A.}
$f: [0,1]^4 \mapsto S$, such that  $X=_d X^*=[X^*_{i,j}],$ where $X^*_{i,j}=f(\alpha, \xi_i,\eta_j, \lambda_{i,j})$.
\end{theorem}

When the array $X$ consists of a single row or a single column, we get a special case, which is equivalent to the Hewitt-Savage theorem (and includes 
de Finetti's theorem):

\begin{theorem}
An infinite $S$-valued sequence  $Z$ is exchangeable if and only if there exists a measurable function 
$f: [0,1]^2 \mapsto S$ and
i.i.d.\ random variables $\alpha, \xi_1,\xi_2,\ldots,$ all uniformly distributed on $[0,1]$, such that 
$Z=_d [f(\alpha,\xi_i)]_{i=1}^\infty$.
\end{theorem}
Note that for any fixed $\alpha=c_0$, the sequence $[f(c_0,\xi_i)]_{i=1}^\infty$ is i.i.d., so with a random $\alpha$ we indeed obtain an i.i.d.\ mix.
Comparing with the formulations  of Theorems \ref{thm1} and \ref{thm2}, observe that here the potentially complicated mixing measure $\mu$ is replaced by the simple random variable $\alpha$, which  is {\em uniform} on $[0,1]$.
Of course, the potential complexity of $\mu$ does not simply ``evaporate," it is just shifted to the function $f$.

\medskip\medskip
\noindent
{\bf de Finetti's theorem for Markov chains.} Diaconis and Freedman \cite{diaconis} created a version of de Finetti's Theorem for Markov chains.
The mixture of  Markov chains can be interpreted similarly to other sequences, as a Markov chain is just a special sequence of random variables.

To elaborate the conditions, consider random variables taking values in a countable state space $I$. Let us call two fixed sequences $a=(a_1,\ldots, a_n)$ and $b=(b_1,\ldots, b_n)$ in $I$ {\em equivalent} if $a_1=b_1$, and the number of $i\rightarrow j$ transitions occurring in $a$ is the same as the number of $i\rightarrow j$ transitions occurring in $b$, for every $i,j\in I$.

Let $X=(X_1,X_2,\ldots)$ be a sequence of random variables over $I$. We say that $X$ is {\em recurrent} if for any starting state $X_1=i$, the sequence 
returns to $i$ infinitely many times, with probability 1. Then the Markov chain version of de Finetti's Theorem, proved by 
Diaconis and Freedman \cite{diaconis}, can be formulated as follows:

\begin{theorem} {\bf (Markov chain version of de Finetti's theorem)}
Let $X=(X_1,X_2,\ldots)$ be a recurrent sequence of random variables over a countable state space $I$. If 
$$
\Pr(X_1=a_1, \ldots,X_n=a_n) \= \Pr(X_1=b_1, \ldots,X_n=b_n)
$$
for any $n$ and for any equivalent sequences $a=(a_1,\ldots, a_n)$, $b=(b_1,\ldots, b_n)$, then $X$ is a mixture of Markov chains.
Furthermore the mixing measure is uniquely determined.
\end{theorem}

\section{The Case of Finite Exchangeable Sequences}
\label{finite}

As already mentioned in Section~\ref{intro}, de Finetti's Theorem does not necessarily hold for finite sequences. There exist, however, 
related results for the finite case, as well. Below we briefly review three fundamental theorems.

\subsection{Approximating a Finite Exchangeable Sequence by an i.i.d.\ Mixture}

Even though  de Finetti's Theorem may fail for finite sequences,
intuition suggests that a finite, but very long sequence will likely behave similarly to an infinite one. This intuition is made precise by a 
result of Diaconis and Freedman \cite{diaconis2}. It provides a sharp bound for the total variation distance between the joint distribution
of exchangeable random variables $X_1,\ldots,X_k$ and the closest mixture of i.i.d.\ random variables. The distance is measured by the 
{\em total variation distance}.  The total variation distance between distributions $P$ and $Q$ is defined as $$d_{TV}(P,Q)=2 \sup_A |P(A)-Q(A)|.$$

\begin{theorem} Let $X_1,\ldots,X_k,X_{k+1},\ldots, X_n$ be an exchangeable sequence of random variables, taking values in an arbitrary 
measurable space $S$. Then the total variation distance between the distribution of $(X_1,\ldots,X_k)$ and of the 
closest mixture of i.i.d.\ random variables is at most $2|S|k/n$ if $S$ is finite, and at most $k(k-1)/n$ if  $S$ is infinite.
\end{theorem}

Observe that the distance bound depends on both $k$ and $n$, and it becomes small only if $k/n$ is small. Thus, if the sequence to be approximated is long 
(i.e., $k$ is large), then this fact in itself does not bring the sequence $(X_1,\ldots,X_k)$ close to an i.i.d.-mix. 
In order to claim such a closeness, we need that $(X_1,\ldots,X_k)$ is {\em extendable} to a significantly longer exchangeable sequence $(X_1,\ldots,X_n)$.

\subsection{Exact Expression of a Finite Exchangeable Sequence by a Signed Mixture}

Another interesting result on the finite case is due to Kerns and Sz\'ekely \cite{kerns}. They proved that any finite
exchangeable sequence, taking values from an arbitrary measurable space, can always be expressed {\em exactly} as an i.i.d.\ mix.
This would not hold in the original setting. But the twist that Kerns and Sz\'ekely have introduced is 
that the mixing measure is a so called {\em signed measure}. The latter means that it may also take negative values.
In the notation recall that ${\cal P}(S)$ denotes the set of all probability distributions on $S$.

\begin{theorem} \label{thm7} 
Let $X_1,\ldots,X_n$ be a sequence of exchangeable random variables, taking values from an arbitrary measurable space $S$. 
Then there exists a signed measure $\nu$ on ${\cal P}(S)$, such that for 
any measurable $A_1,\ldots,A_n \subseteq S$ the following holds:
\begin{equation}\label{signed}
\Pr(X_1\in A_1,\ldots,X_n\in A_n) \= \int_{{\cal P}(S)} \pi(A_1)\cdot\ldots\cdot \pi(A_n)\, d\nu(\pi)
\end{equation}
where $\pi$ runs over ${\cal P}(S)$, integrated according to the signed measure $\nu$. 
\end{theorem}
Here the mixing measure $\nu$ does not have to be unique, in contrast to the traditional versions of the theorem.

A harder question, however, is this: comparing with the traditional versions, the right-hand side of (\ref{signed}) means that  
$\pi$ is drawn according to a signed measure from  ${\cal P}(S)$. What does this mean from the probability interpretation point of view?

Formally, the integral on the right-hand side of (\ref{signed}) it is just a mixture (linear combination, with weights summing to 1) 
of the values $\pi(A_1)\cdot\ldots\cdot \pi(A_n)$, where $\pi$ runs over ${{\cal P}(S)}$. The only deviation from the classical case is 
that some $\pi\in {{\cal P}(S)}$ can be weighted with negative weights.
Thus, {\em formally}, everything is in order, we simply deal with a mixture of {\em probability distributions,} allowing negative weights, but insisting that 
at the end a non-negative function must result. 
However, if we want to interpret it as a mixture of {\em random sequences,} rather than just probability distributions, then the signed measure amounts to a selection via a probability distribution incorporating {\em negative probabilities.} 

What does it mean? How can we pick a value of a random variable with negative probability? To answer this meaningfully is not easy. 
There are some attempts in the literature to interpret negative probabilities, for a short introduction see Sz\'ekely \cite{szekely}. Nevertheless, 
it appears that negative probabilities are neither widely accepted in probability theory, nor usually adopted
in applications, apart from isolated attempts. Therefore, we rather stay with the formal interpretation: ``drawing"
$\pi\in {{\cal P}(S)}$ according to a signed measure for the integral just means taking a mixture (linear combination) of probability distributions 
with weights summing to 1, also allowing negative weights, while insisting that the result is still a non-negative probability distribution.
This makes Theorem~\ref{thm7} formally correct, avoiding  troubles with interpretation. Nevertheless, the interpretation still remains a challenging 
philosophical problem, given that Theorem~\ref{thm7} has been  the only version to date that provides an exact expression of the {\em distribution} any finite 
exchangeable sequence as a mix of i.i.d\ distributions, but it does not correspond to the mixture of random {\em sequences} in the usual (convex) sense.  

\subsection{Exact Finite Representation as a Mixture of Urn Sequences}
\label{mixurn}

Another interesting result about finite exchangeable sequences is that they can be expressed as a mixture (in the usual convex sense) of so-called {\em urn sequences}, explained below. It seems, this provides  the most direct analogy of de Finetti's Theorem for the finite case, yet this result did not receive the attention it deserves, as pointed out by  Carlier, Friesecke, and V\"ogler \cite{carlier}. The idea goes back to de Finetti \cite{definetti3}. Later it was used by several authors at various levels of generality as a proof technique, rather than a target result in itself, see, e.g.,  Kerns and Sz\'ekely \cite{kerns},
so it did not become a ``named" theorem.
Finally, the most general version, which applies to arbitrary random variables, appears in the book of Kallenberg (see \cite{kallenberg},
Proposition 1.8).

{\em Urn sequences} constitute a simple model of generating random sequences. As the most basic version, imagine an urn in which we place $N$ balls, 
and each ball has a certain color. We randomly draw the balls from the urn one by one, and observe the obtained random sequence of colors. 
We can distinguish two basic variants of the process: after drawing a ball, it is put back in the urn (urn process with replacement), or it is not put back
(urn process without replacement).

Consider the following simple example. Let us put $N$ balls in the urn, $K$ black and $N-K$ white balls. If we randomly draw them with replacement, 
then an i.i.d.\ sequence is obtained, in which each entry is black with probability $K/N$, and white with probability $(N-K)/N$.
 The length of the sequence can be arbitrary (even infinite), as the drawing can continue indefinitely. 

On the other hand, if we do this experiment {\em without replacement,} then the maximum length of
the obtained sequence is $N$, since after that we run out of balls. 
The probability that among the first $n\leq N$ draws (without replacement) there are precisely $X$ black balls follows the {\em hypergeometric distribution} 
(see, e.g., Rice \cite{rice}) given by  
\begin{equation}\label{hyper}
\Pr(X=k)\=\frac{{K \choose k}{N-K \choose n-k}}{{N\choose n}}.
\end{equation}
For our purposes the important variant is the case {\em without replacement}, and with $n=N$, that is, all the balls are drawn out of the urn. Then 
the obtained sequence has length $N=n$. Note that it cannot be i.i.d., as it contains precisely $K$ black and $N-K$ white balls. But otherwise it is 
completely random, so the distribution of $X$ is the same as it were in an i.i.d.\ sequence, conditioned on including precisely $K$ black balls. 

The number of colors can be more than two, even infinite. The obtained random sequence is still 
similar to an i.i.d.\ one, with the difference that each color occurs in it a fixed number of times. We can then formulate the general definition of 
the urn sequences of interest to us. For short description, let us first introduce some notations. The set $\{1,\ldots,n\}$ is abbreviated by $[n]$, and 
the family of all permutations of $[n]$ is denoted by $\Sigma_n$. If a permutation $\sigma$ is applied to a sequence $X=(X_1,\ldots,X_n)$, 
then the resulting sequence is denoted by $\sigma(X)$, which is an abbreviation of $(X_{\sigma(1)}, \ldots, X_{\sigma(n)})$.
We also use the following naming convention: 

\noindent {\bf Convention 1 (Uniform random permutation)}
{\em Let $\sigma\in \Sigma_n$ be a permutation. 
We say that $\sigma$ is a  {\em uniform random permutation}, if it is chosen from the uniform distribution over 
$\Sigma_n$.}

Now the urn sequences of interest to us are defined as follows:

\begin{definition}\label{urn} {\bf (Urn sequence)}
Let $x=(x_1,\ldots,x_n)$ be a deterministic sequence, each $x_i$ taking values from a set $S$, and let $\sigma\in \Sigma_n$ be a uniform
random permutation.
Then $X=\sigma(x)$ is called an {\em urn sequence.}
\end{definition}
Here each $x_i$ represents the color of a ball, allowing repeated occurrences. The meaning of $\sigma(x)$ is simply that we list 
the balls in random order. Note that due to the random permutation, we obtain a random sequence, even though $x$ is deterministic. Now we can state the result, after Kallenberg \cite{kallenberg}, but using our own notations:

\begin{theorem}\label{urnthm} {\bf (Urn representation)}
Let $X=(X_1,\ldots,X_n)$ be a finite exchangeable sequence of random variables, each $X_i$ taking values in a
measurable space $S$.
Then $X$ can be represented as a mixture of urn sequences. 
Formally, there exists a probability measure $\mu$ on $S^n$ (mixing measure), 
such that for any $A\subseteq S^n$
$$
\Pr(X\in A) \= \int_{S^n}\Pr\lor \sigma_x(x)\in A \ler d\mu(x)
$$
holds, where $\sigma_x\in \Sigma_n$ is a uniform random permutation, drawn independently for every $x\in S^n$.
\end{theorem}

Observe that Theorem~\ref{urnthm} shows a direct analogy to Theorem~\ref{thm1}, replacing the i.i.d.\ Bernoulli sequence with a finite urn sequence,
giving us the finite length analogy of de Finetti's Theorem.
In the special case when $S=\{0,1\}$, using the hypergeometric distribution formula (\ref{hyper}), we can specialize it to the following result,
resembling the conditional independence form of de Finetti's Theorem, given in Theorem~\ref{thm1b}:

\begin{theorem} \label{thmK}
Let $X_1,X_2,\ldots,X_N$ be a finite sequence of $\{0,1\}$-valued exchangeable random variables. Then there exists a random variable $\eta$,
taking values in $\{0,1,\ldots,N\}$, such that for every $n\in [N]$ and $k\in \{0,1,\ldots,n\}$, the following holds:
\begin{equation}\label{mix2}
\Pr\left(\sum_{i=1}^n X_i=k \;\Largabs\; \eta=K\right) \= \frac{{K \choose k}{N-K \choose n-k}}{{N\choose n}}.
\end{equation}
Furthermore, $\eta$ is given as the number of ones in the sequence, representing the empirical distribution:
$$
\eta \= \sum_{i=1}^N X_i.
$$
\end{theorem}
Theorem~\ref{thmK} says: given that the length-$N$ exchangeable sequence contains $K$ ones, it behaves precisely as an urn sequence that contains $K$ ones.
This also provides a simple algorithm to generate the exchangeable sequence: first pick $\eta$ from its distribution, and whenever $\eta=K$, then 
generate an urn sequence with $K$ ones. The distribution of $\eta$ (the mixing measure) can be obtained as the empirical distribution of ones in the 
original sequence. The sequence generated this way will be statistically indistinguishable from the original.

\section{A Decomposition Theorem for General Finite Sequences}
\label{general}

In all known versions of de Finetti's Theorem, a sequence of rather special properties is represented as a mixture of simpler sequences.
In most cases the target sequence is exchangeable. Although there are some exceptions (some of them are listed in Section~\ref{modifying}), 
the target sequence is {\em always} assumed to satisfy some rather strong symmetry requirement.

Now we raise the question: is it possible to eliminate {\em all } symmetry requirements? That is, can we express an {\em arbitrary} sequence of 
random variables as a mixture of simpler ones? Surprisingly, the answer 
is in the affirmative, with one condition: our method can only handle finite sequences. The reason is that we use uniform 
random permutations, and they do not exist over an infinite sequence. 
On the other hand, we deal with  completely arbitrary random variables, taking values in any measurable space. 

With a general target sequence, the component sequences clearly cannot be restricted to i.i.d., or to urn sequences, since they are all exchangeable,
and the mixture of exchangeable sequences cannot create non-exchangeable ones. Then which class of sequences should the components be taken from? We introduce a class  that we call {\em elementary sequences}, which will do the job. In the definition we use the notation $\alpha\circ \beta$ for the superposition 
(composition) of two permutations, with the meaning $(\alpha\circ \beta)(x)=\alpha(\beta(x))$.

\begin{definition}\label{elementary} {\bf (Elementary sequence)}
Let $x=(x_1,\ldots,x_n)$ be a deterministic sequence, each $x_i$ taking values from a set $S$, and let $\alpha, \beta\in \Sigma_n$ be uniform
random permutations, possibly not independent of each other. Then $X=(\alpha\circ \beta)(x)$ is called an {\em elementary sequence.} 
\end{definition}

Observe the similarity to Definition~\ref{urn}. The only difference is that  in an elementary sequence the permutation is the composition of {\em two }
uniform random permutations,  while in the urn sequence we only use a single uniform random permutation. Of course, if $\alpha$ and $\beta$ 
in Definition~\ref{elementary} are independent of each other, then their superposition would remain a uniform random permutation, giving back
Definition~\ref{urn}.
On the other hand, if they are not independent, then we may get a sequence that is not an urn sequence. 

Let us note that not every sequence is elementary. This follows from the observation that if we fix any $a\in S$, then the number of times $a$ occurs in 
 $(\alpha\circ \beta)(x)$ is constant (which may be 0). The reason is that permutations do not change the number of occurrences of $a$, so its occurrence number 
 remains the  same as in $x$, which is constant. On the other hand, in an arbitrary random sequence this occurrence number is typically random, not constant, so elementary sequences form only a small special subset of  all random sequences. In fact, as we prove later in Lemma~\ref{occur}, the constant occurrence 
 counts actually {\em characterize} elementary sequences. To formalize this, let us introduce the following definition:
 
\begin{definition}\label{occurcount}{\bf (Occurrence count)} Let $X=(X_1,\ldots,X_n)\in S^n$ be a sequence and $a\in S$. Then $F(a,X)$ 
denotes the number of times $a$ occurs in $X$, that is,
$$
F(a,X)\= |\{i\;|\; X_i=a\}|
$$
\end{definition}

The next definition deals with the case when a fixed total ordering $\prec$ is given on $S$. 

\begin{definition}{\bf (Ordered sub-domain, order respecting measure)}\label{order}
The subset of $S^n$ containing all ordered $n$-entry sequences with respect to some total ordering $\prec$ on $S$ is called the {\em ordered sub-domain} 
of $S^n$, denoted by ${\rm Ord}(S^n)$:
$$
{\rm Ord}(S^n) \= \{(x_1,\ldots,x_n)\;|\; x_i\in S, \; (\forall i), \; x_1\prec \ldots \prec x_n\}.
$$
A probability measure on $S^n$ is called {\em order respecting (for the ordering $\prec$),} if $\mu(A)=0$ holds for every measurable set $A\subseteq S^n$, whenever
$A\cap  {\rm Ord}(S^n)=\emptyset$.
\end{definition}

Now we are ready to state and prove our representation theorem for arbitrary finite sequences of random variables.

\begin{theorem}\label{repthm}
Let $X=(X_1,\ldots,X_n)$ be an arbitrary  finite sequence of random variables, each $X_i$ taking values in a
measurable space $S$.
Then $X$ can be represented as a mixture of elementary sequences. 
Formally, there exist a probability measure $\mu$ on $S^n$ (mixing measure),
such that for any measurable $A\subseteq S^n$
\begin{equation}\label{XA}
\Pr(X\in A) \= \int_{S^n} \Pr\lor (\alpha\circ\beta_x)(x)\in A \ler d\mu(x)
\end{equation}
holds, where $\alpha, \beta_x\in \Sigma_n$ are 
uniform random permutations, possibly not independent of each other, and $\beta_x$ is drawn independently for each $x\in S^n$.
 Furthermore, the claim remains true if the mixing measure  $\mu$ is restricted to be order respecting  
for a total ordering $\prec$ on $S$ (see {\rm Definition~\ref{order}}). In that case, the representation is given by the formula
\begin{equation}\label{XA2}
\Pr(X\in A) \= \int_{{\rm Ord}(S^n)} \Pr\lor (\alpha\circ\beta_x)(x)\in A \ler d\mu(x).
\end{equation}
\end{theorem}

For the proof we need two lemmas. The first is a folklore result, stating that  if an arbitrary sequence (deterministic or random, with any distribution) 
is subjected to a uniform random permutation, independent of the sequence, then the sequence becomes exchangeable. We state it below as a lemma for further reference.

\begin{lemma} \label{IRRarb}
Applying an independent uniform random permutation to an arbitrary finite sequence gives an exchangeable sequence.
\end{lemma}

\noindent {\bf Proof.} Let $X=(X_1,\ldots,X_n)$ be an arbitrary finite sequence, taking values from a set $S$. 
Let $Y=(Y_1,\ldots,Y_n)$ be the sequence obtained by applying an independent uniform random permutation $\alpha$ to $X$, i.e., $Y=\alpha(X)$.
Pick $k\leq n$ distinct indices $i_1,\ldots,i_k\in [n]$. Then for any $a_1,\ldots,a_k\in S$ we can write
\begin{equation}\label{atlag}
\Pr(Y_{i_1}=a_1,\ldots,Y_{i_k}=a_k) \= \frac{1}{{n\choose k} k!} \sum_{j_1,\ldots,j_k\in [n],\, \mbox{\footnotesize distinct}} 
\Pr(X_{j_1}=a_1,\ldots,X_{j_k}=a_k).
\end{equation}
The reason is that under the independent uniform random permutation any set of $k$ distinct indices have equal chance to take the place of $i_1,\ldots,i_k$, 
and there are 
${n\choose k} k!$ such sets. As a result, the average obtained on the right-hand side of (\ref{atlag}) does not depend on the specific 
$i_1,\ldots,i_k$ values, only on $k$. Therefore, $\Pr(Y_{i_1}=a_1,\ldots,Y_{i_k}=a_k)$ depends only on $k$, but not on 
$i_1,\ldots,i_k$. This is precisely one of the equivalent definitions of an exchangeable sequence.

\hfill $\spadesuit$

The second lemma expresses the fact that a uniform random permutation can ``swallow" any other permutation, making their composition also
a uniform random permutation. 

\begin{lemma} \label{swallow}
Let $\sigma, \gamma\in \Sigma_n$ be two permutations, such that 
\vspace*{-1mm}
\begin{itemize}\itemsep0.5mm
\item $\sigma$ is a uniform random permutation
\item $\gamma$ is an arbitrary permutation (deterministic or random, possibly 
non-uniform, and possibly dependent on the sequence to which it is applied)
\item $\sigma$ and $\gamma$ are independent of each other. 
\end{itemize}
Then $\sigma\circ\gamma$ is a uniform random permutation.
\end{lemma}

\noindent {\bf Proof.}
Let $\xi=(\xi_1, \ldots,\xi_n)$ be the sequence (deterministic or random) to which the permutation $\gamma$ is applied. Fix an index 
$k\in [n]$, and let $\nu$ be the index to which $\gamma$ maps the index $k$, i.e., 
$\gamma(k)=\nu$. Note that $\nu$ may possibly be random, non-uniform, and dependent on $\xi$. 
Let us express the probability that $\sigma$ maps $\nu$ into a fixed index $i$:
\begin{equation}\label{perm}
\Pr(\sigma(\nu)=i)\= \sum_{\sigma_0}\Pr(\sigma(\nu)=i\;|\;\sigma=\delta_0)\Pr(\sigma=\delta_0),
\end{equation}
where the summation runs over all fixed permutations $\delta_0$ of $[n]$. Observe that $\Pr(\sigma=\delta_0)=1/n!$,
as $\sigma$ is uniform and $\delta_0$ is fixed. Furthermore, 
$$
\Pr(\sigma(\nu)=i\;|\;\sigma=\delta_0) \= \Pr(\delta_0(\gamma(k))=i\;|\;\sigma=\delta_0).
$$
From this, using the independence of $\sigma$ and $\gamma$, we obtain
$$
\Pr(\sigma(\nu)=i\;|\;\sigma=\delta_0) \= \Pr(\delta_0(\gamma(k))=i\;|\;\sigma=\delta_0) \= \Pr(\delta_0(\gamma(k))=i)\= \Pr(\delta_0(\nu)=i).
$$
 Then we can continue (\ref{perm}) as
\begin{equation}\label{perm2}
\Pr(\sigma(\nu)=i)\= \frac{1}{n!} \sum_{\delta_0}\Pr(\delta_0(\nu)=i) \=
\frac{1}{n!} \sum_{j=1}^n \sum_{\delta_0}\Pr(\delta_0(j)=i\;|\;\nu=j)\Pr(\nu=j).
\end{equation}
In the above expression, the event $\{\delta_0(j)=i\}$ involves only fixed values, so it is not random, it happens either with probability 1 or 0,
depending solely on whether $\delta_0(j)=i$ or not. As such, it is independent of the condition $\{\nu=j\}$, so we have 
$\Pr(\delta_0(j)=i\;|\;\nu=j) = \Pr(\delta_0(j)=i)$, whenever the conditional probability is defined, i.e., $\Pr(\nu=j)>0$. If $\Pr(\nu=j)=0$, 
then the conditional probability is undefined, 
but in this case the term cannot contribute to the sum, being multiplied with $\Pr(\nu=j)=0$. Thus, we can continue (\ref{perm2}) as
$$
\Pr(\sigma(\nu)=i)\= \frac{1}{n!} \sum_{j=1}^n \sum_{\delta_0}\Pr(\delta_0(j)=i)\Pr(\nu=j). 
$$
Here the sum $\sum_{\delta_0}\Pr(\delta_0(j)=i)$ is the number of permutations that map a fixed $j$ into a fixed $i$. The number of such permutations
is $(n-1)!$, as the image of $j$ is fixed at $i$, and any permutation is allowed on the rest. This yields
$$
\Pr(\sigma(\nu)=i)\= \frac{1}{n!} \sum_{j=1}^n\underbrace{\sum_{\delta_0}\Pr(\delta_0(j)=i)}_{=(n-1)!}\Pr(\nu=j) \=
\frac{(n-1)!}{n!} \,\underbrace{\sum_{j=1}^n \Pr(\nu=j)}_{=1} \= \frac{1}{n}.
$$
Thus, we obtain $\Pr(\sigma(\nu)=i)\=1/n$, which means that the position to which $\nu=\gamma(k)$ is mapped by $\sigma$ is uniformly distributed over $[n]$,
no matter how $\nu$ was selected, and how it depended on $\xi$. This holds for every $k$, making $\sigma\circ\gamma$ a uniform random permutation.

\hfill $\spadesuit$

\medskip
Before turning to the proof of Theorem~\ref{repthm}, let us point out that a
consequence of the above lemma is interesting on its own right:

\begin{corr}\label{corr1} Any permutation (deterministic or random) can be represented as the composition of two uniform random permutations.
Formally, let $\gamma\in \Sigma_n$ be an arbitrary permutation, deterministic or random; if random, then drawn from an arbitrary distribution. Then there
exist two uniform random permutations $\alpha,\beta\in \Sigma_n$ (possibly not independent of each other), such that $\alpha\circ\beta=\gamma$. 
\end{corr}

\noindent {\bf Proof.} Let $\sigma\in \Sigma_n$ be a uniform random permutation, independent of $\gamma$. Then by Lemma~\ref{swallow}, 
the permutation $\beta=\sigma\circ\gamma$ becomes a uniform random permutation. Set $\alpha=\sigma^{-1}$, which is also a uniform random permutation. 
Further, let $\epsilon$ denote the identity permutation that keeps everything in place.
Then we can write
$$
\alpha\circ\underbrace{\sigma\circ\gamma}_{=\beta} \= \underbrace{\sigma^{-1} \circ\sigma}_{=\epsilon}\circ \gamma \= \gamma,
$$
yielding $\alpha\circ\beta=\gamma$. As $\alpha,\beta$ are both uniform random permutations (possibly not independent of each other), this proves the claim.

\hfill $\spadesuit$

The above corollary also provides an opportunity to characterize elementary sequences:

\begin{lemma}\label{occur}{\bf (Characterization of elementary sequences)}
A sequence $X=(X_1,\ldots,X_n)\in S^n$ is elementary if and only if for any $a\in S$ the occurrence count $F(a,X)$ (see {\rm Definition~\ref{occur}}) is
constant.
\end{lemma}

\noindent {\bf Proof.} If $X$ is elementary, then, by Definition, it can be represented as $(\alpha\circ\beta)(x)$, where $\alpha,\beta\in \Sigma_n$
are uniform random permutations (possibly not independent), and $x\in S^n$ is a deterministic sequence. Since 
no permutation can change occurrence counts, and $F(a,x)$ is constant, due to $x$ 
being deterministic, therefore, $F(a,X)$ remains constant for any $a\in S$. 

Conversely,
assume $F(a,X)$ is constant for any $a\in S$. Let $a_1,\ldots,a_k\in S$ be the distinct elements for which $F(a_i,X)>0, \; i\in [k]$.
Clearly, $k\leq n$, since there can be at most $n$ distinct elements in $X$, and the identity of these elements is fixed,
due to the constant value of $F(a,X)$ for any $a\in S$. Let $y$ be the deterministic sequence that 
contains $a_1,\ldots,a_k$, each one repeated $F(a_i,X)$ times. That is,
$$
y=(\underbrace{a_1,\ldots,a_1}_{F(a_1,X)\, {\rm times}},\underbrace{a_2,\ldots,a_2}_{F(a_2,X)\, {\rm times}}, \ldots,
\underbrace{a_k,\ldots,a_k}_{F(a_k,X)\, {\rm times}})
$$
Then we have $F(a,y)=F(a,X)$ for every $a\in S$. Thus, $X$ and $y$ contain the same elements, with the same multiplicities, 
just possibly in different order. That is, $X$ is a permutation of $y$, possibly a random permutation, which may depend on $y$. 
Let $\gamma_y\in \Sigma_n$ be the permutation that implements $X=\gamma_y(y)$. Then 
by Corollary~\ref{corr1}, the permutation $\gamma_y$ can be represented as $\alpha\circ\beta=\gamma_y$, where $\alpha,\beta\in \Sigma_n$ are
uniform random permutations, possibly not independent of each other, and they may also depend on $y$. But no matter
what dependencies exist, Corollary~\ref{corr1} provides that 
$X=\gamma_y(y)$ can be represented as $X=(\alpha\circ\beta)(y)$ for some uniform random permutations $\alpha,\beta\in \Sigma_n$
and a deterministic sequence $y$, proving that $X$ is indeed elementary.

\hfill $\spadesuit$

\medskip
\noindent {\bf Proof of Theorem~\ref{repthm}.} Let us apply a uniform random permutation $\rho\in \Sigma_n$ to $X$, such that $\rho$ and $X$ 
are independent. This results in  a new sequence 
$Y=\rho(X)$. By Lemma~\ref{IRRarb}, the obtained $Y$ is an exchangeable sequence. Then by Theorem~\ref{urnthm} we have that  
$Y$ can be represented as a mixture of urn sequences. That is, there exists
a probability measure $\mu$ on $S^n$, 
such that for any $A\subseteq S^n$
\begin{equation}\label{YA2}
\Pr(X\in A) \= \int_{S^n}\Pr\lor \sigma_x(x)\in A \ler d\mu(x)
\end{equation}
holds, where $\sigma_x\in \Sigma_n$ is a uniform random permutation, drawn independently for every $x\in S^n$.
 This representation means that $Y$ can be produced by drawing $x$ from the the 
mixing measure $\mu$, and drawing a uniform random permutation $\sigma_x\in \Sigma_n$, and  then outputting $\sigma(x)$. 

Now, instead of outputting $\sigma_x(x)$, let us first permute it by $\rho^{-1}$. Thus, we  output $(\rho^{-1}\circ \sigma_x)(x)$. 
Observe that if $\rho$ is a uniform random permutation, then so is  $\rho^{-1}$, which we denote by $\alpha$. This makes
the resulting $(\rho^{-1}\circ \sigma_x)(x)=(\alpha\circ \sigma_x)(x)$ an elementary sequence. 
Applying $\alpha$ to the mixture  means that each component sequence $\sigma_x(x)$ is permuted 
by  $\alpha$. But then the result is also permuted by $\alpha$, since it does not matter whether the components are permuted first, 
and then one of them is selected, or the selection is made first and the result is permuted afterwards with the same permutation.

Applying $\alpha$ in the above way, we obtain the sequence $\alpha(Y)=\rho^{-1}(Y)$ as the result. Thus we can re-write (\ref{YA2}) as 
\begin{equation}\label{YA3}
\Pr\lor \alpha(Y)\in A \ler \= \int_{S^n} \Pr\lor (\alpha\circ \sigma_x)(x)\in A \ler d\mu(x).
\end{equation}
Now we observe that $\alpha(Y) = \rho^{-1}(\rho(X)) = X$.  Then we can continue (\ref{YA3}) as 
\begin{equation}\label{XA1}
\Pr\lor X\in A \ler \= \int_{S^n} \Pr\lor (\alpha\circ \sigma_x)(x)\in A \ler d\mu(x),
\end{equation}
which is precisely the formula (\ref{XA}) we wanted to prove, just using the notation $\sigma_x$ instead of $\beta_x$.

Consider now the case when $\mu$ is order respecting  for some ordering $\prec$ on $S$.
Let $\gamma_x\in \Sigma_n$ be the permutation that orders $x=(x_1,\ldots,x_n)\in S^n$ according to $\prec$, that is 
$\gamma_x(x)=(x^*_1,\ldots,x^*_n)$ where $x^*_1\prec\ldots\prec x^*_n$ is the ordered version of $x_1,\ldots,x_n$. Let $\delta\in \Sigma_n$
be a uniform random permutation, chosen independently of $\gamma_x$. Then  $\delta$ and $\gamma_x$ satisfy the conditions of Lemma~\ref{swallow}. Therefore, by Lemma~\ref{swallow}, 
$\delta\circ\gamma_x$ is a uniform random permutation. Introducing the notation $\beta_x=\delta\circ\gamma_x$,  we obtain 
from the already proven formula (\ref{XA})
$$
\Pr\lor X\in A \ler \= \int_{S^n} \Pr\lor (\alpha\circ \beta_x)(x)\in A \ler d\mu(x),
$$
where $\alpha,\beta_x$ are uniform random permutations, and $\beta_x$ is chosen independently for
each $x\in S^n$. Since $\mu$ is order respecting (see Definition~\ref{order}), it is enough to restrict the integration to 
the set ${\rm Ord}(S^n)$, giving us the formula (\ref{XA2}). This completes the proof.

\hfill $\spadesuit$

\section{Application of de Finetti Style Theorems in Random Network Analysis}
\label{applic}

Large, random networks, such as wireless ad hoc networks, are often described by various types of random graphs, primarily by 
geometric random graphs. A frequently used model is when each node 
of the network is represented as a random point in some planar domain, and two such nodes are connected by an edge (a network link) if they are within a given
distance from each other. This basic model has many variants: various domains may occur, different probability distributions of the node positions within the domain may be used, a variety of distance metrics is possible, etc. Note that it falls in the category of static random graph models, which is our focus here, 
in contrast to evolving ones (for a survey of random graph models, see e.g., Drobyshevskiy and Turdakov \cite{droby}). Let us now formalize 
what we mean by a general random graph model.
 
\begin{definition} {\bf (Random graph models)} \label{RGM}
Let $X=(X_1,X_2,X_3,\ldots)$ be an infinite sequence of random variables, each taking its values from a fixed domain $S$, 
which is an arbitrary measurable space.
A {\em random graph model over $S$} is a function $\bf G$ that maps $X$ into a sequence of graphs:
$${\bf G}(X) \= (G_1,G_2,G_3\ldots).$$
If $X$ is restricted to a subset $C\subseteq S^\infty$, then we talk about a {\em conditional random graph model}, denoted by 
${\bf G}(X\,|\,  C)$.
\end{definition}
Note that even though the random graph model ${\bf G}(X)$ depends on the infinite sequence $X$, the individual graphs $G_n$ typically depend
only on an initial segment of $X$, such as $(X_1,\ldots,X_n)$.

Regarding the condition $C$, a very simple variant is when $C=C_1\times C_2\times C_3\times\ldots$, where $C_i\subseteq S$, and we independently restrict each $X_i$ to fall into $C_i$. Note, however, that $C$ may be much more complicated, possibly not reducible to individual restrictions on each $X_i$. 

A most frequently occurring case is when the points (the components of $X$)  are selected from the same distribution independently, that is, they are i.i.d. The reason is that allowing dependencies makes the analysis too messy. To this end, let us define i.i.d.-based random graph models:

\begin{definition} {\bf (i.i.d.-based random graph models)} \label{iidRGM}
If the entries of $X$ in {\rm Definition~\ref{RGM}} are i.i.d. random variables, 
then we call ${\bf G}(X)$ an {\em i.i.d.-based random graph model over $S$}, and ${\bf G}(X\,|\,  C)$ is called an 
{\em i.i.d.-based conditional random graph model over $S$}.
\end{definition}

The most commonly used and analyzed static random graphs are easily seen to fall in the category of 
i.i.d.-based random graph models. Typical examples are Erd\H os-R\'enyi random graphs (when each edge is added independently with some
probability $p$), different variants of geometric random graphs, random intersection graphs,
and many others. On the other hand, sometimes the application provides natural reasons for considering
dependent points, as shown by the following example.

\medskip
\noindent {\sc Example 3.} Consider a wireless ad hoc network. 
Let each node be a point drawn independently and uniformly from  the unit square. Specify a transmission radius $r>0$, and 
connect two nodes whenever they are within distance $r$ (Note: $r=r_n$ may depend on the number of nodes.) 
But allow only those systems of points for which the arising graph has diameter (in terms of  graph distance) of at most some value $D=D_n$, 
which may again depend on $n$. The conditioning makes the points dependent. Nevertheless, the restriction is reasonable if we want the
network to experience limited delays in end-to-end transmissions.

This example (and many possible similar ones) shows that there can be good reasons  to deviate from the standard i.i.d.\ assumption.
On the other hand, most of the analysis results build on the i.i.d.\ assumption. How can we bridge this gap? Below we show an approach that 
is grounded in de Finetti style theorems, and provides a tool that can come handy in the analysis of conditional random graph models. 

\begin{theorem}\label{CRGMthm}
Let $S$ be a Borel measurable space, and let $\cal P$ be a property of random graph models. Fix an i.i.d.-based random graph model $\bf G$ over $S$, and 
assume that ${\bf G}(X)$ has property ${\cal P}$, regardless of the value of $X$, with probability 1. 
Let $C$ represent a condition with $\Pr(X\subseteq C)>0$.
Then ${\bf G}(X\,|\, C)$  also has property ${\cal P}$,  regardless of the value of $X$, with probability 1.
\end{theorem}

\noindent {\bf Proof.} Let $Y$ be a random variable that has the conditional distribution of $X$, given $C$. That is,
for every measurable set $A$
$$
\Pr(Y\subseteq A)\= \Pr(X\subseteq A\;|\; X\subseteq C).
$$
Note that $Y$ may not remain i.i.d. However, we show that $Y$ is still exchangeable. Let $\sigma$ be any permutation. Then we can write
\begin{equation}\label{sigmaY}
\Pr(\sigma(Y)\subseteq A)\= \Pr(\sigma(X)\subseteq A\;|\; \sigma(X)\subseteq C) \= \frac{\Pr(\sigma(X)\subseteq A\cap C)}{\Pr(\sigma(X)\subseteq C)}.
\end{equation}
Since $X$ is i.i.d., therefore, $X=_d\sigma(X)$, i.e., they have the same distribution, making them statistically indistinguishable. This implies that for any measurable set $B$
the equality $\Pr(\sigma(X)\subseteq B) = \Pr(X\subseteq B)$ holds. Using it in (\ref{sigmaY}), we get that 
$$
\Pr(\sigma(Y)\subseteq A) \= \frac{\Pr(\sigma(X)\subseteq A\cap C)}{\Pr(\sigma(X)\subseteq C)} \=
\frac{\Pr(X\subseteq A\cap C)}{\Pr(X\subseteq C)} \= \Pr(Y\subseteq A)
$$ for any permutation $\sigma$, which means that $Y$ is exchangeable. Here we also used that $\Pr(X\subseteq C)>0$, so the denominator 
does not become 0.

Recall now that by the Hewitt-Savage Theorem, 
an infinite sequence $Y$ of $S$-valued exchangeable random variables is an $S$-valued i.i.d.-mix, whenever  
$S$ is a Borel measurable space.
For each component $X$ of this i.i.d.\ mix 
we can take the function $\bf G$ to obtain a random graph model ${\bf G}(X)$. After taking the mixture, this results in ${\bf G}(Y)$. The reason is that
applying the function to each sequence $X$ first and then selecting one of  them must yield the same result as first selecting one of them (by the same 
mixing measure), and applying the function to the selected sequence $Y$.

As a result of the above reasoning, we get ${\bf G}(X)=_d{\bf G}(Y)$, i.e., the two random graph models have the same distribution. 
But $Y$ was chosen such that it has the conditional distribution of $X$, given $C$. Therefore, we have ${\bf G}(Y)=_d{\bf G}(X\,|\, C)$.
Thus, we obtain ${\bf G}(X\,|\, C)=_d{\bf G}(X)$. Since, by assumption, ${\bf G}(X)$ has property ${\cal P}$,
regardless of the value of $X$, with probability 1,
therefore, ${\bf G}(X\,|\, C)$  also has property ${\cal P}$, regardless of the value of $X$, with probability 1. 
The reason for we need that $\cal P$ does not depend on $X$ (with probability 1) is that when we mix various realizations of 
$X$, they should all come with the same property $\cal P$, otherwise a mixture of properties would result.
This completes the proof.

\hfill $\spadesuit$

The above result may sound very abstract, so let us illustrate it by two examples. 

\medskip
\noindent
{\sc Example 4.} It follows from the results of Farag\'o \cite{farago_peva} that every i.i.d.-based geometric random graph has the following property: 
\begin{quote} 
If the graph is asymptotically connected (that is, the probability of being connected approaches 1 as the number of nodes tends to infinity), 
then the average degree must tend to infinity. 
\end{quote}
Let us choose this as property $\cal P$. 
One may ask: does this property remain valid in conditional models over the same geometric domain? 
We may want to know this in more sophisticated  models, such as the one presented in Example 3. Observe that the above property $\cal P$ 
satisfies the condition that it holds regardless of the value of $X$ (with probability 1), where $X$ represents the random points on which the 
geometric random graph model ${\bf G}(X)$ is built. 
Therefore, by Theorem~\ref{CRGMthm}, the property remains valid for ${\bf G}(X\,|\, C)$, as well,
no matter how tricky and complicated condition is introduced, as long as the condition holds with positive probability (even when $n\rightarrow\infty$). 
Note that this cuts through a lot of complexity that may otherwise arise if we want to prove the same claim directly
from the specifics of the model.

\medskip
\noindent
{\sc Example 5.} Consider the variant of Erd\H os-R\'enyi random graphs, where each edge is added independently with some
probability $p$. These random graphs are often denoted by $G_{n,p}$, where $n$ is the number of vertices. For constant $p$, they fit in our general random graph 
model concept, choosing $X$ now as an i.i.d.\ sequence of Bernoulli random variables, representing the edge indicators. Let $\kappa(G_{n,p}),
\lambda(G_{n,p})$ and $\delta(G_{n,p})$ denote the vertex connectivity, edge connectivity and minimum degree of $G_{n,p}$, respectively.
All these graph parameters become random variables in a random graph. A nice (and quite non-trivial) result from the theory of random graphs (see Bollob\'as \cite{bollobas}) is that for any $p$, the following holds:
\begin{equation}\label{boll}
\lim_{n\rightarrow\infty} \Pr \lor \kappa(G_{n,p})\= \lambda(G_{n,p}) \= \delta(G_{n,p}) \ler \= 1.
\end{equation}
The intuitive meaning of (\ref{boll}) is that asymptotically both types of connectivity parameters are determined solely by the minimum degree.
The minimum degree always provides a trivial lower bound both for $\kappa(G_{n,p})$ and  $\lambda(G_{n,p})$, and in a random graph asymptotically 
they both indeed hit this lower bound, with probability 1.

Now we may ask: what happens if we introduce some condition? Let $\cal G$ be a subset of graphs that represents a condition that 
$G_{n,p}$ satisfies with some constant probability 
$q,\; 0<q<1$, for every $n$. That is $\Pr(G_{n,p} \in {\cal G})=q$, for every $n$. Observe that (\ref{boll}) holds regardless of the value of $X$,
because (\ref{boll}) is valid for every $p$. Therefore, it can be used as property $\cal P$ in Theorem~\ref{CRGMthm}.
Thus, if we condition on $G_{n,p}$ falling in $\cal G$, the relationship (\ref{boll}) still 
remains true, by  Theorem~\ref{CRGMthm}. 

Note that if $\cal G$ is complicated, it may be very hard to prove it directly from the model that (\ref{boll}) remains true
under the condition $G_{n,p} \in {\cal G}$. Fortunately, our result cuts through this complexity. It is also interesting to note that in this case
$X$ is a Bernoulli sequence, so for this case it would be enough to use the original de Finetti Theorem in the proof, rather than the more 
powerful Hewitt-Savage Theorem.

\section{Conclusion}

The first part of the paper  reviews some results regarding the decomposition of random sequences into a mixture of simpler ones. 
This line of research started with the classic theorem of de Finetti, and later it was extended and generalized in numerous directions. 
Since it is not considered very well known in the Engineering/Computer Science community, we provide more details than what is usual 
in the introductory parts of articles. Then we have presented a new representation theorem in Section~\ref{general}, which covers cases not considered before. 
Finally, in Section~\ref{applic}, we have demonstrated that de Finetti style results can provide unexpected help in the analysis of random networks.


\medskip\medskip
\noindent
{\Large\bf Appendix A: Measurable Spaces and Related Concepts}

\medskip
\noindent
For the sake of self-containment, we briefly explain/summarize some background concepts that are referred to in the paper.

\smallskip
\noindent {\bf Measurable space.}
A {\em measurable space} is a pair ${\cal S}=(S,{\cal A})$, where $S$ is a set,
and $\cal A$ is special set system, a {\em $\sigma$-algebra} over $S$. 
The $\sigma$-algebra $\cal A$ is a family of subsets of $X$, with the following properties: it contains $X$ 
itself, it is closed under taking complements and countable unions. (These properties imply that $\cal A$ also contains the empty set and is closed under countable intersections, as well.) If $S$ is finite or countably infinite, then $\cal A$ often simply contains all subsets of $S$, but this is not necessary. The sets that are contained 
in $\cal A$ are called {\em measurable sets}. 

Why do we need to distinguish measurable sets? Because in some situations, typically for 
non-countable models, we cannot avoid dealing with non-measurable sets (see later, under the heading {\em Non-measurable sets}).
For simplicity, when we talk about a  measurable space, we  often just denote it by its underlying set $S$, 
rather than the more precise ${\cal S}=(S,{\cal A})$ notation; this usually does not cause any confusion.

\smallskip
\noindent {\bf Measurable function.}
If ${\cal S}_1=(S_1,{\cal A}_1)$ and ${\cal S}_2=(S_2,{\cal A}_2)$ are two (not necessarily different) measurable spaces, then a function
$f:S_1\mapsto S_2$ is called a {\em measurable function} if for any $A\in {\cal A}_2$, the set of elements that the function maps into $A$ constitute 
a measurable set in ${\cal S}_1$. That is, 
$\{x\,|\, f(x)\in A\}\in {\cal A}_1$.
The set $f^{-1}(A)=\{x\,|\, f(x)\in A\}$ is called the {\em pre-image} of $A$. Informally, the condition is often expressed this way: the function is measurable if and only if the pre-image of any measurable set is also measurable. 

\smallskip
\noindent {\bf Isomorphic measurable spaces.}
Two measurable spaces ${\cal S}_1=(S_1,{\cal A}_1)$ and ${\cal S}_2=(S_2,{\cal A}_2)$ are called {\em isomorphic}, if there exists a bijection (1--1 onto function) $f:S_1\mapsto S_2$, such that both $f$ and its inverse are measurable. 

\smallskip
\noindent {\bf Borel measurable space.}
The {\em Borel subsets} of {\bf R} (the set of real numbers) are the sets that arise by repeatedly taking countable unions, countable intersections 
and relative complements (set differences) of open sets. A measurable space is called a {\em Borel measurable space} if it is isomorphic to 
a Borel subset of {\bf R}.

\smallskip
\noindent {\bf Measure space.}
Note that no measure is included in the definition of a measurable space. If a measure is also added, then it becomes a {\em measure space,} 
not to be confused with a {\em measurable} space. A {\em measure} is a function that assigns a real number to every measurable set, such that 
certain axioms are satisfied. Specifically, if ${\cal S}=(S,{\cal A})$ is a measurable space, then a function $\mu:{\cal A}\mapsto {\bf R}$ is a measure,
if it is non-negative, $\mu(\emptyset)=0$, and is countably additive. The latter means that for every countable collection of sets $A_1,A_2,\ldots\in \cal A$
it holds that $\mu (\cup_{i=1}^\infty A_i) = \sum _{i=1}^\infty \mu(A_i)$. Then the triple ${\cal M}=(S,{\cal A}, \mu)$ is referred to as a measure space.

\smallskip
\noindent {\bf Probability space.}
A {\em probability measure} is a measure with the additional requirement that the measure of the whole space is 1. If this is satisfied, then the 
arising measure space 
is referred to as a {\em probability space} or a {\em probability triple}. The parts are often denoted differently from the notation 
$(S,{\cal A}, \mu)$ of a general measure space. A frequently used notation for a probability space is $(\Omega,{\cal F}, P)$, where $\Omega$ is the 
set of possible outcomes (elementary events), $\cal F$ is the collection of events, and $P$ is the probability measure.

\smallskip
\noindent {\bf Non-measurable sets.}
The subsets that belong to the $\sigma$-algebra $\cal F$ of subsets in a probability space represent the possible events we want to deal with. 
Why do we bother with a $\sigma$-algebra rather than simply allowing all subsets as  possible events? We can certainly do it if $\Omega$ is finite. 
In the infinite case, however, we need to be careful. For example, if $\Omega={\bf R}$, then provably there exists no measure on {\em all} subsets that
satisfies the axioms of a probability measure; there are always non-measurable sets, even though they tend to be contrived (the proof requires the Axiom of Choice). 
 
 \smallskip
\noindent {\bf Random variable.}
A {\em random variable}, in the most general setting, is a measurable function from a probability space to a measurable space.  Let us illustrate it with an example. 
Let $\Omega$ be the set of all infinite bit sequences, containing infinitely many 1-bits. Each such infinite string is a possible outcome of an experiment. 
For each such string let us assign a real number in $[0,1]$, which is obtained by viewing the bit string as the binary expansion of the number,
after a leading 0 and the decimal point
(it will be a 1--1 mapping, due to requiring infinitely many 1-bits). 
For a set of strings, let the probability measure of the set be some standard measure of the size of the corresponding set of real numbers,
such as the Lebesgue measure. Then the $\sigma$-algebra $\cal F$ of events is the family of those string sets that map into Lebesgue measurable subsets of 
$[0,1]$. (This does not contain all subsets, as there are non-measurable sets, albeit contrived ones.) To define a random variable, let us chose the set of all non-negative integers as the target measurable space, allowing all subsets in its $\sigma$-algebra. 
Let a random variable defined by the function that maps a bit string into the integer that tells
how many 1-bits are among the first 100 bits of the string. It is not hard to see that this satisfies the general definition of a measurable function. Therefore,
it indeed correctly defines a random variable. 


\end{document}